\documentclass{amsart}
\usepackage{amssymb}
\usepackage{graphicx}
\usepackage{url}
 \newtheorem{theorem}{Theorem}
 \newtheorem{lemma}[theorem]{Lemma}
 \newtheorem{proposition}[theorem]{Proposition}
 \newtheorem{corollary}[theorem]{Corollary}
 \newtheorem{remark}[theorem]{Remark}
 \newtheorem{example}[theorem]{Example}
 \newtheorem{definition}[theorem]{Definition}
 \newtheorem{conjecture}[theorem]{Conjecture}

 \newtheorem{question}[theorem]{Question}

\newcommand{\bpr}{\begin{proof}}
\newcommand{\epr}{\end{proof}}
\newcommand{\beq}{\begin{equation}}

\newcommand{\bThm}{\begin{theorem}}
\newcommand{\eThm}{\end{theorem}}
\newcommand{\blem}{\begin{lemma}}
\newcommand{\elem}{\end{lemma}}
\newcommand{\bpro}{\begin{proposition}}
\newcommand{\epro}{\end{proposition}}
\newcommand{\bcor}{\begin{corollary}}
\newcommand{\ecor}{\end{corollary}}
\newcommand{\brem}{\begin{remark}}
\newcommand{\erem}{\end{remark}}
\newcommand{\bexa}{\begin{example}}
\newcommand{\eexa}{\end{example}}
\newcommand{\bdf}{\begin{definition}}
\newcommand{\edf}{\end{definition}}
\newcommand{\bcon}{\begin{conjecture}}
\newcommand{\econ}{\end{conjecture}}
\newcommand{\bque}{\begin{question}}
\newcommand{\eque}{\end{question}}

\newcommand{\comment}[1]{}

\title{On the strength of invariants from stacks of virtual knots}
\author{Blake K. Winter}
\address{Mathematics Dept., University of Findlay, Findlay, OH}
\address{bw3073@findlay.edu}

\keywords{knots, virtual knots, algebraic invariants}
\subjclass{57K12}

\date{}

\begin{document}
\thispagestyle{empty}

\begin{abstract}
For any virtual link, a class of new links can be defined called stacks, in which copies of the virtual link are placed on top of one another. The resulting virtual link depends only on the virtual isotopy class of the original link, and the fundamental group of such a link may be used to detect whether the link is nontrivial and whether it is nonclassical in some cases. We show that the groups constructed using this method are sufficient to distinguish all the Kishino knots from the unknot and from one another, as well as calculating their Jones polynomials. 
\end{abstract}

\maketitle

%%%%%%%%%%%%%%%%%%%%%%%%%%%%%%%%%%%%%%%%%%%%%%%%%%%
%
\section{Introduction}

Virtual links were introduced by Kauffman \cite{Kauff}, and given a geometric interpretation by Kuperberg \cite{Kup}. The concept of a stack for a virtual link was introduced in \cite{BW1}. The general idea is to consider the virtual link as a link in a thickened surface, and then to glue copies of the link in the thickened surface below the original link, possibly with some copies having the orientation of the thickening reversed. The virtual isotopy class of the new link depends only on the virtual isotopy class of the original link, and thus the new link allows for constructing invariants for the original.

In particular, the group and quandle of the stack form an invariant for the original link. In \cite{BW1} it is shown that this can be nontrivial even if the group of the original link is trivial. However, it is remarked that this is not true for two layer stacks of a certain Kishino knot (herein named $K_{switch}$) which is not distinguished from the unknot by the groups of these two layer stacks. The motivation for introducing stacks in \cite{BW1} was to attempt to construct topological homotopy invariants that would be more powerful than the usual group and quandle, and as such the question was posed therein whether it was possible to use homotopy invariants to distinguish $K_{switch}$ from the unknot by using more layers in the stack. Our main result is that this is possible using a three layer stack. To our knowledge this is the first use of a homotopy invariant to detect the nontriviality of this knot, although this is possible using biquandle counting invariants \cite{Nelson1} or colored Jones polynomials. In addition the Jones polynomials of the stacks are calculated, which distinguish all seven cases using two layer stacks.

\section{The Stack Invariant}

Let $L$ be a virtual link. Then following Kuperberg \cite{Kup} we may regard $L$ as a link (also designated by $L$) in a thickened surface $F\times I$. Virtual isotopy is then the equivalence relation generated by two types of modifications. We may (smoothly) isotope $L$ in $F\times I$. We also regard $L_1 \subset F_1 \times I$ as equivalent to $L_2 \subset F_2 \times I$ if there is a (smooth) map $f:F_1 \rightarrow F_2$ such that $f\times id_I : F_1\times I \rightarrow F_2 \times I$ is an embedding and $f\times id_I(L_1)=L_2$.

Let $L$ be a virtual link represented by $L\subset F\times I$. We define the vertical reflection $L_v$ of $L$ to be the link obtained by reversing the orientation of the thickening line segment, or equivalently, if $L:\prod_1 ^n S^1 \rightarrow F\times I$ is written as $L(x) = (y, t)$, then $L_v (x) = (y, 1-t)$, where $(y, t)$ is a point in $F\times I$. This has the effect of switching each crossing to interchange which strand is the overcrossing strand and which is the undercrossing strand.

For a sequence $\langle a_i\rangle .$ of $+$ and $-$ symbols, of length $m$, we define the stack invariant of $L$ to be the link $S_{a_i}$ to be the link constructed as follows. Let $L_i\subset F_i\times I$ be a copy of $L$ if $a_i = +$ and a copy of $L_v$ if $a_i=-$. Then identify $F_i\times \{0 \}$ with $F_{i+1}\times \{1 \}$. Finally reparametrize the resulting thickened surface so that the thickening is identified with the unit interval in the obvious way. We refer to $L_i$ as the $i$-th layer of the stack, and refer to $m$ as the number of layers.

In the case $S_{+-}(L)$ \cite{BW1} refers to this particular stack as the \emph{vertical double} of $L$, written $VD(L)$, and we will continue that notation. We will also write simply $DD(L)$ (the \emph{double}) for $S_{++}(L)$.

\section{Kishino Knot Calculations}

There are seven knots which we will consider here, which exhaust the Kishino knots up to horizontal reflection and orientation. Since horizontal reflection and orientation will not affect the link groups of the stacks (or the nontriviality of the Jones polynomials considered later), this is sufficient for our purposes, namely, to test the strength of the homotopy invariants of stacks by using the Kishino knots. Note that the vertical reflections may change the link groups of the stacks. However, also note that vertical reflections have related stacks: if $L$ and $L'$ are vertical reflections of one another, then $S_{-+}(L)$ is identical to $S_{+-}(L')$, for example. Nonetheless, we consider the vertical reflections of the Kishino knots separately, since we are interested in whether invariants from stacks can distinguish them.

In particular, we will consider $K$, $K_{switch}$, $K_{alt}$, $K_v$, $K_5$, $K_6$, and $K_7$. The groups are calculated by checking crossings, and the presentations simplified using Magma \cite{Magma}, which was also used to count inequivalent group epimorphisms. Note that all these knots have cyclic knot groups except $K_v$, which has a group given by

\begin{flalign*}
   &\langle a,
    d|
    d  a d^{-1} a d^{-1} a^{-1}  = 1\rangle . 
		\end{flalign*}
Note that this alone shows that not only $K_v$ but also $K$ is nontrivial and both are distinct from the other five: $K$ and $K_v$ are vertical reflections of one another, but none of the others have nontrivial groups for themselves or their vertical reflections.
\begin{figure}
		\centering
			\includegraphics[scale=0.5]{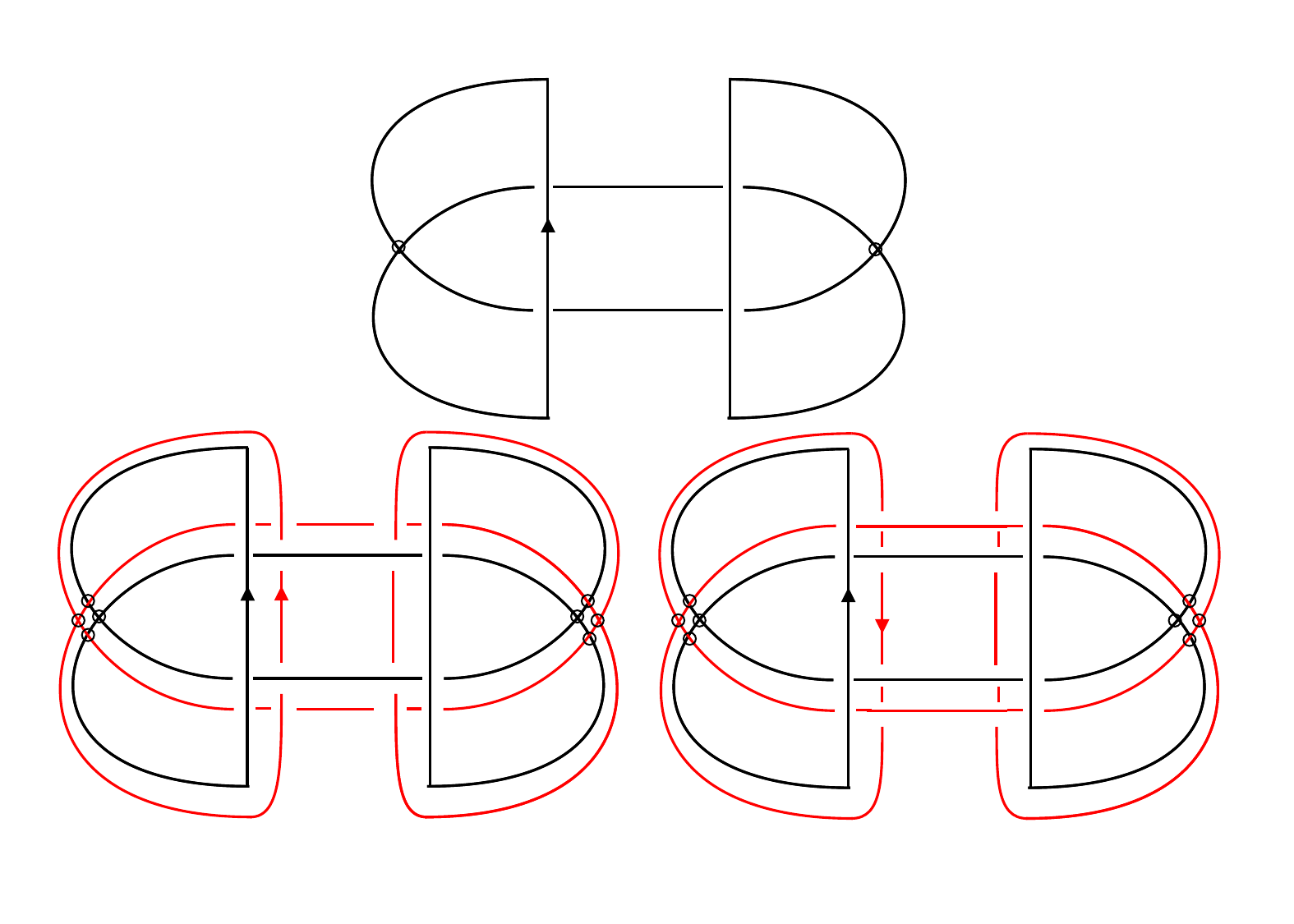}
		\caption{The knot $K$ with the vertical double on the right and the double on the left.}
		\label{K1}
\end{figure}
\begin{figure}
		\centering
			\includegraphics[scale=0.5]{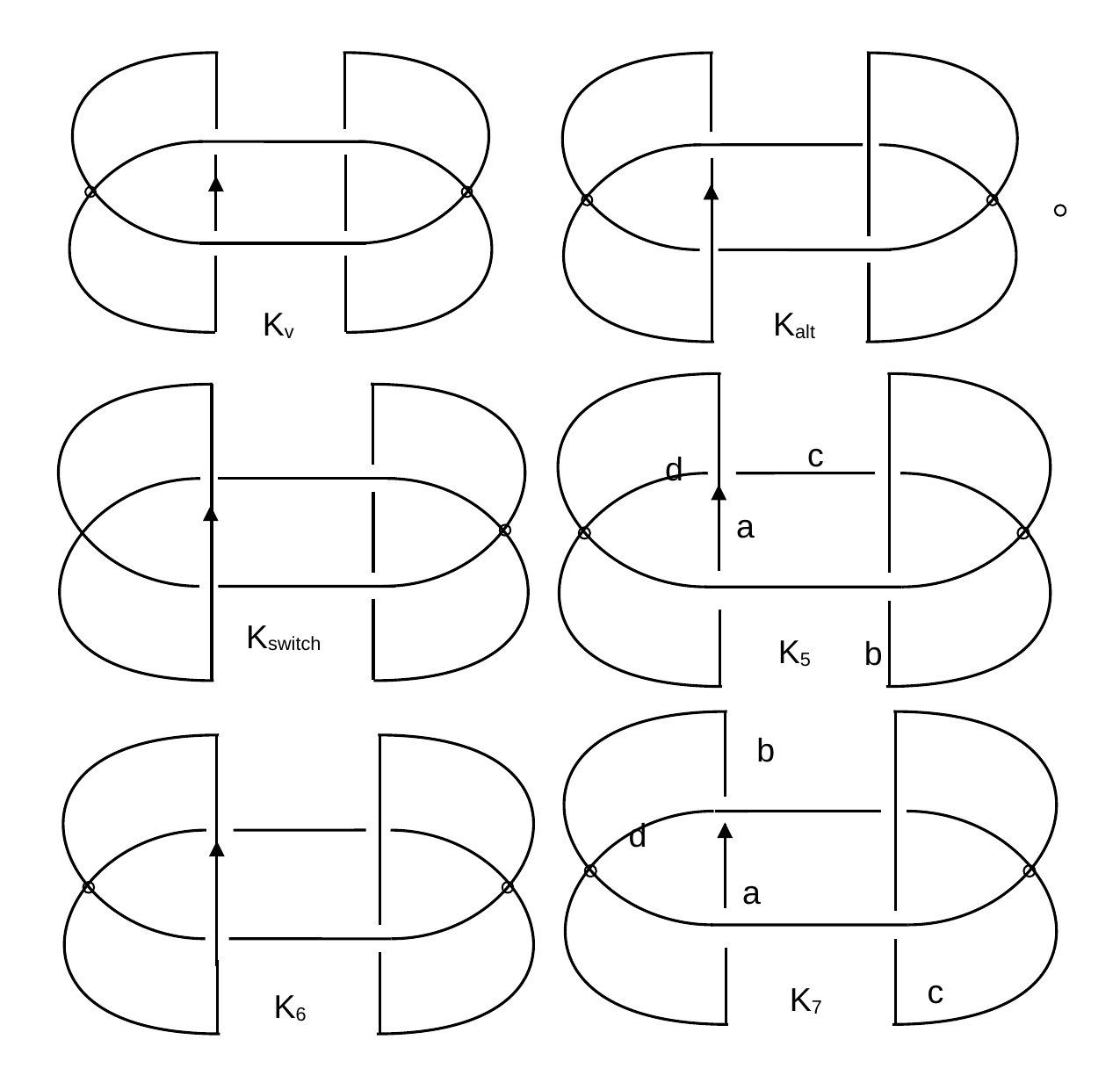}
		\caption{The other six Kishino knots considered in this paper.}
		\label{Kset}
\end{figure}
%\begin{figure}
%		\centering
%			\includegraphics[scale=0.5]{Kswitch.pdf}
%		\caption{The knot $K_{switch}$ with the vertical double on the right and the double on the left.}
%		\label{Kswitch}
%\end{figure}
%\begin{figure}
%		\centering
%			\includegraphics[scale=0.5]{Kalt.pdf}
%		\caption{The knot $K_{alt}$ with the vertical double on the right and the double on the left.}
%		\label{Kalt}
%\end{figure}
%\begin{figure}
%		\centering
%			\includegraphics[scale=0.5]{Kv.pdf}
%		\caption{The knot $K_{v}$ with the vertical double on the right and the double on the left.}
%		\label{Kv}
%\end{figure}
%\begin{figure}
%		\centering
%			\includegraphics[scale=0.5]{K5.pdf}
%		\caption{The knot $K_{5}$ with the vertical double on the right and the double on the left.}
%		\label{Kv}
%\end{figure}
%\begin{figure}
%		\centering
%			\includegraphics[scale=0.5]{K6.pdf}
%		\caption{The knot $K_{6}$ with the vertical double on the right and the double on the left.}
%		\label{Kv}
%\end{figure}
%\begin{figure}
%		\centering
%			\includegraphics[scale=0.5]{K7.pdf}
%		\caption{The knot $K_{7}$ with the vertical double on the right and the double on the left.}
%		\label{Kv}
%\end{figure}

From the diagrams we can calculate the link groups:

For $S_{+-}(K)$:
%\begin{flalign*}
%&\langle a, b, c, d, e, f, g, h, i, j, k, l, m |\\
%&c=b^a, c=d^k, k=l^a,\\
%&e=d^a, m=b^a, m=l^e,\\
%&i=h^a, e=f^a, i=j^e,\\
%&g=f^k, k=j^a, g=h^a\rangle . 	
%\end{flalign*}

%This simplifies to
\begin{flalign*}
   &\langle a,
    e,
    k |
    e a^{-1}  e^{-1} k^{-1} a^{-1} e k^{-1} a = 1\rangle . 	
		\end{flalign*}

This has 75 inequivalent epimorphisms to $S_5$.

For $S_{++}(K)$:
%\begin{flalign*}
%	&\langle a, b, c, d, e, f, g, h, i, j, k, l, m |\\
%	&c=b^a, d=c^a, e=d^m, 
%	\\&e=f^i, f=g^a, g=h^a, 
%	\\&i=h^a, j=i^a, k=j^g, \\
%	&k=l^c, l=m^a, m=b^a\rangle .
%\end{flalign*}
%This simplifies to
\begin{flalign*}
   &\langle a,
    c,
    g |
    g c^{-1} a^{-1} c a c g^{-1} a^{-1} g^{-1} a = 1\rangle . 	
		\end{flalign*}

This has 75 inequivalent epimorphisms to $S_5$.

%5or $S_{+-}(K_{switch})$:
%\begin{flalign*}
%&\langle a, b, c, d, e, f, g, h, i, j, k, l, m |\\
%&c=b^a, c=d^k, k=l^a,\\
% &m=b^a, m=l^e, e=d^a,\\
%&f=e^i, f=g^a, i=h^a,\\
%&j=i^a, g=h^a, j=k^g\rangle .
%\end{flalign*}
%For $S_{++}(K_{switch})$:
%\begin{flalign*}
%&\langle a, b, c, d, e, f, g, h, i, j, k, l, m |\\
%&c=b^a, l=m^a, k=l^c,\\
%&m=b^a, d=c^a, e=d^m,\\
%&j=i^e, i=h^a, e=f^a,\\
%&k=j^a, g=h^a, f=g^k\rangle .
%\end{flalign*}

The groups for $S_{+-}(K_{switch})$ and $S_{++}(K_{switch})$  both simplify to free groups on two generators.

For $S_{+-}(K_{alt})$:
%\begin{flalign*}
%&\langle a, b, c, d, e, f, g, h, i, j, k, l, m |\\
%&c=b^a, k=l^a, k=j^c,\\
%&d=c^a, m=b^a, m=l^d,\\
%&e=d^i, e=f^a, i=h^a,\\
%&j=i^a, g=h^a, g=f^j\rangle .
%\end{flalign*}

%This simplifies to
\begin{flalign*}
&\langle a, c |
    a c a^{-2} c a c a^{-1} c^{-1} a^2
    c^{-1} a^{-1} c a^{-2} c a c^{-1} a^{-1}
    c^{-1} a^2 c^{-1} a^2 c a^{-2} c a c^{-1}\\
    &a^{-1} c^{-1} a^2 c^{-1} a^{-2} c^{-1} a^2 c
    a^{-2} c a c a^{-1} c^{-1} a^2 c^{-1}
    a^{-2} c a^{-2} c a c a^{-1} c^{-1} a^2
    c^{-1} = 1 \rangle .
\end{flalign*}

This has 12 inequivalent epimorphisms to $S_5$.

For $S_{++}(K_{alt})$:
%\begin{flalign*}
%&\langle a, b, c, d, e, f, g, h, i, j, k, l, m |\\
%&c=b^a, d=c^l, l=m^a,\\
%&e=d^a, f=e^m, m=b^a,\\
%&f=g^a, i=h^a, j=i^f,\\
%&l=k^g, k=j^a, g=h^a\rangle .
%\end{flalign*}

%This simplifies to
\begin{flalign*}
&\langle a, c |
a c^{-1} a^{-2} c^{-1} a c a^{-1} c a^2
    c a^{-2} c a^2 c^{-1} a^{-2} c^{-1} a
    c^{-1} a^{-1} c a^2 c a^{-2} c^{-1} a^{-2}\\
    &c^{-1} a c^{-1} a^{-1} c a^2 c a c^{-1}
    a^{-2} c^{-1} a c^{-1} a^{-1} c a^2 c
    a^{-1} c^{-1} a^{-2} c^{-1} a c a^{-1} c
    a^2 c a = 1\rangle .
\end{flalign*}

This has 12 inequivalent epimorphisms to $S_5$.

For $S_{+-}(K_{v})$:
%\begin{flalign*}
%&\langle A, B, C, D, b, c, d, e, f, g, h, i, j, k, l, m |\\
%&B=A^D, D=A^B, D=C^B, B=C^D,\\
%&c=b^a, l=k^c, l=m^a,\\
%&d=c^a, d=e^m, m=b^a,\\
%&i=h^a, f=e^i, f=g^a,\\
%&j=i^a, j=k^g, g=h^a\rangle .
%\end{flalign*}
%This simplifies to
\begin{flalign*}
&\langle B, D, c, g |\\
&D B^{-1} D B^{-1} D^{-1} B = 1,\\
&B^{-2} D c D^{-1} B^2 c^{-1} g B^{-2} D
g^{-1} D^{-1} B^2 g^{-1} c = 1,\\
&B^{-1} D g^{-1} c B^{-1} D^{-1} B c D
B^{-1} D c^{-1} g D^{-1} B D^{-1} g^{-1} D
    = 1\rangle .
\end{flalign*}

This has 241 inequivalent epimorphisms to $S_5$.

For $S_{++}(K_{v})$:
%\begin{flalign*}
%&\langle A, B, C, D, b, c, d, e, f, g, h, i, j, k, l, m |\\
%&B=A^D, D=A^B, D=C^B, B=C^D,\\
%&d=c^k, k=l^a, c=b^a,\\
%&e=d^a, l=m^e, m=b^a,\\
%&i=h^a, e=f^a, j=i^e,\\
%&k=j^a, f=g^k, g=h^a\rangle .                     
%\end{flalign*}
	
%Which simplifies to
\begin{flalign*}
&\langle B, D, e, k |\\
    &B D^{-1} B D^{-1} B^{-1} D = 1,\\
    &B e B^{-1} k^{-1} e B k^{-1} B^{-1} e^{-1}
    k = 1,\\
    &D e D^{-1} k^{-1} D^{-1} B k B e^{-1}
    B^{-1} k^{-1} B^{-1} D k = 1\rangle .                     
\end{flalign*}	

This has 87 inequivalent epimorphisms to $S_5$.

%For $S_{+-}(K_{5})$:
%\begin{flalign*}
%&\langle a, b, c, d, e, f, g, h, i, j, k, l, m |\\
%&c=b^a, c=d^l, l=m^a,\\
%&m=b^a, e=d^a, e=f^m,\\
%&k=j^a, g=f^k, g=h^a,\\
%&i=j^a, i=h^l, l=k^a\rangle .                     
%\end{flalign*}

%For $S_{++}(K_{5})$:
%\begin{flalign*}
%&\langle a, b, c, d, e, f, g, h, i, j, k, l, m |\\
%&c=b^a, j=k^c, k=l^a,\\
%&m=b^a, l=m^d, d=c^a,\\
%&g=f^a, h=g^d, d=e^a,\\
%&j=i^e, e=f^a, i=h^a\rangle .                     
%\end{flalign*}

The groups of $S_{+-}(K_{5})$ and $S_{++}(K_{5})$ both simplify to free groups on two generators.

For $S_{+-}(K_{6})$:
%\begin{flalign*}
%&\langle a, b, c, d, e, f, g, h, i, j, k, l, m |\\
%&c=b^a, c=d^k, k=l^a,\\
%&m=b^a, e=d^a, m=l^e,\\
%&f=e^j, j=i^a, f=g^a,\\
%&k=j^a, h=g^k, h=i^a\rangle .                     
%\end{flalign*}

%This simplifies to

\begin{flalign*}
&\langle a, h |\\
    &h a^{-2} h a h a^{-1} h^{-1} a^2 h^{-1}
    a^{-2} h a h a^{-2} h a h^{-1} a^{-1}\\
    &h^{-1} a^2 h^{-1} a^{-2} h a h a^{-1}
    h^{-1} a^2 h^{-1} a^{-1} h^{-1} a^2 = 1\rangle .                     
\end{flalign*}	

This has 14 inequivalent epimorphisms to $S_5$.

For $S_{++}(K_{6})$:
%\begin{flalign*}
%&\langle a, b, c, d, e, f, g, h, i, j, k, l, m |\\
%&c=b^a, k=l^c, l=m^a,\\
%&m=b^a, e=d^m, d=c^a,\\
%&i=h^e, e=f^a, h=g^a,\\
%&k=j^f, f=g^a, j=i^a\rangle .                     
%\end{flalign*}

%This simplifies to

\begin{flalign*}
&\langle a, c |\\
    &a c^{-1} a^{-1} c a c a^{-1} c^{-1} a^{-1}
    c a c a c^{-1} a^{-1} c^{-1} a c \\
    &a^{-2} c^{-1} a^{-1} c^{-1} a c a c^{-1} 
    a^{-1} c^{-1} a c a^{-1} c^{-1} a^{-1} c a
    c a = 1\rangle .                     
\end{flalign*}	

This has 15 inequivalent epimorphisms to $S_5$.

For $S_{+-}(K_{7})$:
%\begin{flalign*}
%&\langle a, b, c, d, e, f, g, h, i, j, k, l, m |\\
%&c=b^a, l=k^c, l=m^a,\\
%&d=c^a, d=e^m, m=b^a,\\
%&f=e^j, j=i^a, f=g^a,\\
%&h=g^k, h=i^a, k=j^a\rangle .                     
%\end{flalign*}

%This simplifies to

\begin{flalign*}
&\langle a, c |\\
    &a c a^{-1} c^{-1} a c^{-1} a^{-1} c a^{-1}
    c a c^{-1} a c a^{-1} c a c^{-1}\\
    &a^{-2} c a^{-1} c a c^{-1} a c a^{-1}
    c^{-1} a c^{-1} a^{-1} c a^{-1} c^{-1} a
    c^{-1} a = 1\rangle .                     
\end{flalign*}	

This has 15 inequivalent epimorphisms to $S_5$.

For $S_{++}(K_{7})$:
%\begin{flalign*}
%&\langle a, b, c, d, e, f, g, h, i, j, k, l, m |\\
%&c=b^a, d=c^k, k=l^a,\\
%&m=b^a, e=d^a, l=m^e,\\
%&i=h^e, h=g^a, e=f^a,\\
%&j=i^a, k=j^f, f=g^a\rangle .                     
%\end{flalign*}

%This simplifies to

\begin{flalign*}
&\langle a, d |\\
    &d^{-1} a^{-2} d^{-1} a d a^{-1} d a^2 d
    a^{-2} d^{-1} a d^{-1} a^{-2} d^{-1} a d\\
    &a^{-1} d a^2 d^{-1} a^{-2} d^{-1} a d^{-1}
    a^{-1} d a^2 d a^{-1} d a^2 = 1\rangle .                     
\end{flalign*}	

This has 14 inequivalent epimorphisms to $S_5$.

For all the cases where we have reduced the presentation to a free group with one relator, results on one relator groups show that these cannot be written as the free product of the cyclic group on $a$ with another group (see \cite{Group}). Hence they are distinct from the link group of a 2-component unlink. Furthermore, by counting inequivalent epimorphisms to $S_5$, it may be seen that the nontrivial groups for the vertical doubles and the doubles are all distinct from one another (and are also distinct from the free group on two generators).

\begin{table}[]
\begin{tabular}{ll}
$S_{+-}(K)$          & 75   \\
$S_{++}(K)$          & 75   \\
$S_{+-}(K_{switch})$ & free \\
$S_{++}(K_{switch})$ & free \\
$S_{+-}(K_{alt})$    & 12   \\
$S_{++}(K_{alt})$    & 12   \\
$S_{+-}(K_{v})$      & 241  \\
$S_{++}(K_{v})$      & 87   \\
$S_{+-}(K_{5})$      & free \\
$S_{++}(K_{5})$      & free \\
$S_{+-}(K_{6})$      & 14   \\
$S_{++}(K_{6})$      & 15   \\
$S_{+-}(K_{7})$      & 15   \\
$S_{++}(K_{7})$      & 14  
\end{tabular}
\caption{A table showing the inequivalent epimorphisms from the groups to $S_5$. Groups marked "free" are free on two generators.}
\end{table}

For the case of $K_{switch}$ and $K_5$, we can go to a three layer stack to distinguish these from the unknot (and from one another). This positively answers a question in \cite{BW1} about whether $K_{switch}$ can be detected using homotopy groups from stacks.

The group for $S_{+++}(K_{switch})$:

%\begin{flalign*}
%&\langle a, b, c, d, e, f, g, h, i, j, k, l, m, A, B, C, D, E, F, G, H, I, J, K, L, M, N, O, P, Q, R, S, T |\\
%&c=b^a, l=m^a, k=l^c,
%m=b^a, d=c^a, e=d^m,
%j=i^e, i=h^a, e=f^a,
%k=j^a, g=h^a, f=g^k,\\
%&B=A^a, C=B^k, P=Q^C, Q=R^c, R=S^a,
%D=C^a, E=D^m, F=E^S, S=T^e, T=A^a,\\
%&N=M^F, F=G^j, G=H^a, M=L^e, L=K^a,
%O=N^a, P=O^f, H=I^P, I=j^k, J=K^a\rangle .                     
%\end{flalign*}

%This simplifies to

\begin{flalign*}
&\langle a, c, C |\\
    &a^2 c^{-1} a^{-1} c a c a^{-1} C^{-1}
    c^{-1} a^{-1} C a c C a c^{-1} a^{-1}
    c^{-1} a c a^{-2} C^{-1} c^{-1} a^{-1} C^{-1}\\
    &a c C c^{-1} a^{-1} c^{-1} a c a
    c^{-1} a^{-1} c^{-1} a c a^{-1} c^{-1} a^{-1}
    c a c C^{-1} c^{-1} a^{-1} C a c C = 1\rangle .                     
\end{flalign*}	

This has 94 inequivalent epimorphisms to $S_5$.

The group for $S_{+++}(K_5)$ has 928 inequivalent epimorphisms to $S_5$.

$K_1$, $K_{switch}$, and $K_v$ all have trivial Jones polynomials, and require the use of colored Jones to detect, although the other four can be detected by standard Jones polynomials.

$K_{alt}$: $t^{-5} - 2t^{-4} - 2t^{-7/2} + t^{-3} + 2t^{-5/2} + t^{-2}$.

$K_5$: $3 - t^{-3/2} - t^{-1} + t^{-1/2} + t^{1/2} - t - t^{3/2}$.

$K_6$: $-t^{-5/2} + t^{-3/2} + t^{-1}$.

$K_7$: $-t^{-5/2} + t^{-3/2} + t^{-1}$.

Calculating the Jones polynomial of $S_{++}(K_{switch})$ yields $- (1 + t) (1 - 3 t + 3 t^2 - t^3 + 3 t^4 - 3 t^5 + t^6)/t^{7/2}$ \cite{AS}. The PD code was found and sent to Adam Sikora \cite{AS}, who used it to calculate this polynomial in Mathematica using the package KnotTheory$\ \grave{ }$ \cite{KT}. The Jones polynomial of 
$S_{+-}(K_{switch})$ was identical to that of the two-component unlink, however.

For $S_{+-}(K_5)$, the Jones polynomial is $(-1 + t^2 (5 + t - 6 t^2 - 6 t^3 + t^4 + 5 t^5 - t^7))/t^{9/2}$, while for $S_{++}(K_5)$, the Jones polynomial is $(3 + 8 t - 6 t^3 - 6 t^4 + 8 t^6 - 3 t^7)/t^{7/2}$. These were calculated in Mathematica. 

For completeness we list the Jones polynomials of the other five examples considered.

For $S_{++}(K)$, $(-1 + 2 t - 2 t^3 - 2 t^4 + 2 t^6 - t^7)/t^{7/2}$.

For $S_{+-}(K)$, $(-1-t)/t^{1/2}$.
For $S_{++}(K_{v})$, $(-1 + 2 t - 2 t^3 - 2 t^4 + 2 t^6 - t^7)/t^{7/2}$.

For $S_{+-}(K_{v})$, $(-4 + 8 t - 5 t^3 - 5 t^4 + 8 t^6 - 4 t^7)/t^{7/2}$.
For $S_{++}(K_{alt})$, $(-1 + 7 t - 18 t^2 + 20 t^3 - 15 t^4 + 21 t^5 - 19 t^6 + 7 t^7 - 
 4 t^8)/t^{21/2}$.

For $S_{+-}(K_{alt})$, $(-1 + 5 t^2 + t^3 - 6 t^4 - 6 t^5 + t^6 + 5 t^7 - t^9)/t^{9/2}$.

For $S_{++}(K_6)$, $(-3 + 9 t - 7 t^2 + 4 t^3 - 9 t^4 + 6 t^5 - 2 t^6)/t^{11/2}$.
For $S_{+-}(K_6)$, $(1 - 2 t^2 - 2 t^3 + t^5)/t^{5/2}$.

For $S_{++}(K_7)$, $(-3 + 9 t - 7 t^2 + 4 t^3 - 9 t^4 + 6 t^5 - 2 t^6)/t^{11/2}$.

For $S_{+-}(K_7)$, $(-2 + 5 t - 4 t^3 - 4 t^4 + 5 t^6 - 2 t^7)/t^{7/2}$.

We see therefore that the seven Kishino knots are distinguished from each other, and from the unknot, by the Jones polynomials of the two layer stacks $S_{+-}$ and $S_{++}$.

\begin{question}
Determine the relationship, if any, between the (colored) Jones polynomials of a virtual link and the (colored) Jones polynomials of the stacks built from the link.
\end{question}

Note that $K_6$ and $K_7$ have the same Jones polynomial.

\begin{question}
Is it coincidental that for both $K_{switch}$ and $K_5$, both the double and vertical double have trivial group invariants, whereas for the other five, both the double and vertical double have nontrivial group invariants? Does there exist a virtual link $L$ such that the group of the double (or vertical double) is indistinguishable from that of the unlink, but the group of the vertical double (or double) is distinguishing?
\end{question}

The particular interest of this study lies in the fact that, while all of the knots studied here can be classified by means such as colored Jones polynomials or biquandle counting invariants, this is the first time to our knowledge that a classical algebraic topology invariant has been used to classify them. This leads us to a last question:

\begin{question}
Are the groups of stacks as strong (or stronger) invariants as the biquandle or the (colored) Jones polynomial invariants?
\end{question}

%
%%%%%%%%%%%%%%%%%%%%%%%%%%%%%%%%%%%

%%%%%%%%%%%%%%%%%%%%%%%%%%%%%%%%%

\end{document}